\documentclass[fleqn,11pt,a4paper,leqno]{article}

\usepackage{graphicx}
\usepackage{graphics}
\usepackage{latexsym}
\usepackage{amssymb,amsmath}
\usepackage{multicol}
\usepackage{bussproofs}
\usepackage{enumerate}
\usepackage{xcolor}
\usepackage[latin1]{inputenc}
\usepackage{verbatim, color}
\usepackage{enumerate}

\title{proof of associative of $\land$}

\author{}

\date{}

\numberwithin{equation}{section}

\newtheorem{definicion}{Definition}[section]

\newtheorem{definition}[definicion]{Definition}

\newtheorem{Theorem}[definicion]{Theorem}
\newtheorem{theorem}[definicion]{Theorem}
\newtheorem{Corollary}[definicion]{Corollary}

\newtheorem{lema}[definicion]{Lemma}

\newenvironment{Proof}{\noindent\bf Proof \rm}{$\hfill
\square$}

\begin{document}

 \newcounter{thlistctr}
 \newenvironment{thlist}{\
 \begin{list}%
 {\alph{thlistctr}}%
 {\setlength{\labelwidth}{2ex}%
 \setlength{\labelsep}{1ex}%
 \setlength{\leftmargin}{6ex}%
 \renewcommand{\makelabel}[1]{\makebox[\labelwidth][r]{\rm (##1)}}%
 \usecounter{thlistctr}}}%
 {\end{list}}

	\title{Implication Zroupoids and Birkhoff Systems}
	
	\author{Juan M. CORNEJO  and Hanamantagouda P. SANKAPPANAVAR\footnote{The authors wish to dedicate this work to children and their families who fight against cancer.}\ \ }



\maketitle

\begin{abstract}


An algebra $\mathbf A = \langle A, \to, 0 \rangle$, where $\to$ is binary and $0$ is a constant, is called an 
\emph{implication zroupoid} ($\mathcal{I}$-zroupoid, for short) if $\mathbf A$ satisfies the identities:
 $(x \to y) \to z \approx [(z' \to x) \to (y \to z)']'$, where $x' : = x \to 0$,  and  
 $ 0'' \approx 0$.  These algebras generalize De Morgan algebras and $\lor$-semilattices with zero.   Let $\mathcal{I}$ denote the variety of implication zroupoids. 
 For details on the motivation leading to these algebras, we refer the reader to \cite{sankappanavarMorgan2012} (or the relevant papers mentioned at the end of this paper).    
The investigations into the structure of the lattice of subvarieties of $\mathcal{I}$, begun in \cite{sankappanavarMorgan2012}, have continued in \cite{cornejo2016order, cornejo2016semisimple, cornejo2015implication, cornejo2016derived, cornejo2016Bol-Moufang, associativetype,  
cornejo2016weakassociative} and \cite{Implsem}.  
The present paper is a sequel to this series of papers and is devoted to 
making further contributions to the theory of implication zroupoids.

The identity (BR): $x \land (x \lor y) \approx x \lor (x \land y)$ is called the Birkhoff's identity.  The main purpose of this paper is to prove that if $\mathbf A$ is an algebra in  
the variety $\mathcal{I}$, then the derived algebra  $\mathbf A_{mj} := \langle A; \wedge, \vee \rangle$, where $a \land b := (a \to b')'$ and $a \lor b := (a' \land b')'$, satisfies the Birkhoff's identity.   
As a consequence, we characterize the implication zroupoids $\mathbf A$ whose derived algebras  $\mathbf A_{mj}$ 
are Birkhoff systems.
It also follows from the main result that there are bisemigroups that are not bisemilattices but satisfy the Birkhoff's identity, which suggests a more general notion, than Birkhoff systems, of ``Birkhoff bisemigroups'' as bisemigroups satisfying the Birkhoff's identity.  The paper concludes with an open problem on Birkhoff bisemigroups. 
\end{abstract}

\indent identity, Birkhoff system, lattice of subvarieties\\

{\em 2010 AMS subject class:} $Primary: 06D30, 06E75$; $Secondary: 08B15, 20N02,
   03G10$

\section{Introduction}
\

  An algebra $\mathbf A = \langle A, \to, 0 \rangle$, where $\to$ is binary and $0$ is a constant, is called an 
\emph{implication zroupoid} ($\mathcal{I}$-zroupoid, for short) if $\mathbf A$ satisfies the following identities:
\begin{itemize}
\item[{\rm(1)}] \quad $(x \to y) \to z \approx [(z' \to x) \to (y \to z)']'$, where $x' : = x \to 0$,  and  
\item[{\rm(2)}] \quad $ 0'' \approx 0$.   Let $\mathcal{I}$ denote the variety of implication zroupoids.  
\end{itemize}
These algebras generalize De Morgan algebras and $\lor$-semilattices with zero.  For more details on the motivation leading to these algebras, we refer the reader to \cite{sankappanavarMorgan2012} (or the relevant papers mentioned at the end of this paper).  
 
 
The investigations into the structure of the lattice of subvarieties of $\mathcal{I}$, begun in \cite{sankappanavarMorgan2012}, have continued in \cite{cornejo2016order, cornejo2016semisimple, cornejo2015implication, cornejo2016derived, cornejo2016Bol-Moufang, associativetype,  
cornejo2016weakassociative} and \cite{Implsem}.  (It should be noted that in \cite{cornejo2015implication} implication zroupoids were referred to as ``implicator groupoids''.)  
The present paper is a sequel to this series of papers and is devoted to 
making further contributions to the theory of implication zroupoids. 

Throughout this paper we use the following definitions:
\[
	\textup{(M)} \quad 	x \land y := (x \to y')' \quad\text{ and }\quad 
	\textup{(J)} \quad  x \lor y := (x' \land y')'.
\]
With each  $\mathbf A \in \mathcal{I}$, we associate the following algebras:
\[
	\mathbf{A{^{mj}}} := \langle A, \land, \lor, 0 \rangle \quad\text{ and }\quad 
	\mathbf{A{_{mj}}} := \langle A, \land, \lor \rangle.  
\]


\begin{Theorem} {\rm \cite[Corollary 4.6]{cornejo2016derived}}  \label{Theo_Assoc}
	If $\mathbf A \in \mathcal I$ then $\langle A, \wedge \rangle$ and $\langle A, \vee \rangle$ are semigroups.  Hence, $\mathbf{A{_{mj}}}$ is a bisemigroup.
\end{Theorem}

Two of the important subvarieties of $\mathcal{I}$ are:
 $\mathcal{I}_{2,0}$ and $\mathcal{MC}$ which are defined relative to $\mathcal{I}$, respectively, by the following identities: 
\begin{equation} \label{eq_I20}  \tag{{\rm I}$_{2,0}$}
x'' \approx x.
\end{equation}
\begin{equation} \label{eq_MC} \tag{MC}
x \wedge y \approx y \wedge x.
\end{equation}

\begin{definition}
Members of the variety $\mathcal{I}_{2,0}$ are called \emph {involutive}, and members of   
$\mathcal{MC}$ are called \emph {meet-commutative}.  An algebra $\mathbf{A} \in  \mathcal{I}$ is {\rm symmetric} if $\mathbf{A}$ is both involutive and meet-commutative. 
\end{definition}
 Let $\mathcal{S}$ denote the variety of symmetric $\mathcal{I}$-zroupoids.  Thus, $\mathcal{S} = \mathcal{I}_{2,0} \cap \mathcal{MC}$.
The identity
\[
{\rm(BR)} \quad x \land (x \lor y) \approx x \lor (x \land y).
\]
is called the \emph {Birkhoff's identity}.  
This identity, a weakened form of the absorption identities, was introduced by Birkhoff in 1948. In fact,  
 Birkhoff asked in \cite[Problem 7]{birkhoff1948} for an investigation of algebras satisfying the lattice identities without absorption identities but with the identity (BR), which led to the following notion:     
\begin{definition}
A Birkhoff system is a bisemilattice satisfying the Birkhoff's identity {\rm(BR)}.
\end{definition}

Indeed, in response to Birkhoff's problem, there have been a series of papers in the literature revealing the structure of the lattice of subvarieties of the variety of Birkhoff systems; for example, see \cite{harding2012proyective}, \cite{harding2016aprojective}, \cite{harding2016bprojective} and the references therein.  More recently, it was proved in \cite[Theorem 7.3]{cornejo2015implication} that if $\mathbf{A} \in \mathcal{S}$, then 
$\mathbf{A}^{mj}$  is a distributive Birkhoff system, from which it immediately follows that $\mathbf{A}_{mj}$  is a distributive Birkhoff system--a result which will be strengthened in this paper.   


The main purpose of this paper is to prove the following result:
\begin{Theorem}\label{main} If $\mathbf{A}$ is an implication zroupoid, then the bisemigroup $\mathbf{A}_{mj}$ satisfies the Birkhoff identity.
\end{Theorem}
As a consequence, we characterize the implication zroupoids $\mathbf A$ for which  $\mathbf A_{mj}$ 
is a Birkhoff system.

It also follows from the main result, Theorem \ref{main}, that there are bisemigroups that are not bisemilattices but satisfy the Birkhoff's identity, which naturally suggests a more general notion, than Birkhoff systems, of ``Birkhoff bisemigroups'' as bisemigroups satisfying the Birkhoff identity.  This notion seems to be new.  In this new terminology, we can now recast our main theorem as: If $\mathbf{A} \in \mathcal I$, 
then $\mathbf{A}_{mj}$ is a Birkhoff bisemigroup.

The paper concludes with an open problem on Birkhoff bisemigroups. 






\section{Preliminaries}   
 
 In this section we present some preliminary results that will be useful later.

\begin{lema} {\rm \cite[Theorem 8.15]{sankappanavarMorgan2012}} \label{general_properties_equiv}
	The following identities are equivalent in the variety $\mathcal I$:
	\begin{enumerate}[{\rm (a)}]
		\item $0' \to x \approx x$, \label{TXX} 
		\item $x'' \approx x$,
		\item $(x \to x')' \approx x$, \label{reflexivity}
		\item $x' \to x \approx x$. \label{LeftImplicationwithtilde}
	\end{enumerate}
\end{lema}

The following theorem is proved in \cite[Theorem 7.3]{cornejo2015implication}.  
\begin{theorem}\label{theorem_I20_MC}
 Let $\mathbf A \in \mathcal{S}$.   
	Then $\mathbf A_{mj}$   
	satisfies:
	\begin{enumerate}[{ \rm(a)}]
		\item $x \land x \approx x$, \label{010415_01}
		\item $x \lor x \approx x$,
		\item $x \lor y \approx y \lor x$,
		\item $x \land (x \lor y) \approx x \lor (x \land y)$. \label{010415_03} 
	\end{enumerate}
\end{theorem}

\begin{lema}\label{general_properties} \label{general_properties2}
	Let $\mathbf A \in \mathcal I_{2,0}$. Then $\mathbf A$ satisfies:
	\begin{enumerate}[{\rm (1)}]
		\item $x' \to 0' \approx 0 \to x$, \label{cuasiConmutativeOfImplic2}
		\item $0 \to x' \approx x \to 0'$, \label{cuasiConmutativeOfImplic}
		\item $(x \to 0') \to (y \to z) \approx ((0 \to x) \to y) \to z$, \label{260615_01}
		\item $(0 \to x) \to (0 \to y) \approx x \to (0 \to y)$, \label{311014_06}
		\item $0 \to (x \to y) \approx x \to (0 \to y)$, \label{071114_04}
		\item $0 \to (x' \to y)' \approx x \to (0 \to y')$, \label{071114_01}
		\item $(x \to y) \to (y \to z) \approx (0 \to x') \to (y \to z)$, \label{250615_04}
		\item $((x \to y) \to z) \to (z \to u) \approx (0 \to x) \to ((y \to z) \to (z \to u))$, \label{250615_06} 
		\item $x \to y \approx x \to (x \to y)$, \label{031114_04} 
		\item $(y \to x) \to y \approx (0 \to x) \to y$, \label{291014_10}
		\item $(x \to y)' \to (0 \to x)' \approx y' \to x'$, \label{271114_02} 
		\item $(x \to y)' \to y \approx x \to y$, \label{250615_01} 
		\item $[x \to (y \to x)']' \approx (x \to y) \to x$, \label{291014_09}
		\item $x \to ((0 \to x) \to y) \approx x \to y$, \label{250615_22} 
		\item $x \to (y \to x') \approx y \to x'$, \label{281114_01}
		\item $(x \to y) \to y' \approx y \to (x \to y)'$, \label{071114_05}
		\item $((x \to y) \to (z \to x)) \to u \approx (y \to 0') \to ((z \to x) \to u)$, \label{260615_05} 
		\item $(z \to x) \to (y \to z) \approx (0 \to x) \to (y \to z)$, \label{080415_01}
		\item $(x \to y')' \to z \approx x \to (y \to z)$, \label{140715_20}
		\item $0 \to (x \to y')' \approx 0 \to (x' \to y)$, \label{191114_05}
		\item $0 \to (0 \to x)' \approx 0 \to x'$, \label{031114_07}
		\item $[x' \to (0 \to y)]' \approx (0 \to x) \to (0 \to y)'$. \label{031114_06}
	\end{enumerate}
\end{lema}

\begin{Proof}
Items (\ref{cuasiConmutativeOfImplic2}) and (\ref{cuasiConmutativeOfImplic}) are proved in \cite{sankappanavarMorgan2012}.
The proofs of items (\ref{311014_06}), (\ref{071114_04}), (\ref{071114_01}), (\ref{031114_04}), (\ref{291014_10}), (\ref{271114_02}), (\ref{291014_09}), (\ref{281114_01}), (\ref{071114_05}), (\ref{191114_05}), (\ref{031114_07}), (\ref{031114_06}) can be found in \cite{cornejo2016order}.
Items (\ref{260615_01}), (\ref{250615_04}), (\ref{250615_06}), (\ref{250615_01}), (\ref{250615_22}), (\ref{260615_05}), (\ref{080415_01}), (\ref{140715_20}) are proved in \cite{cornejo2016derived}.
\end{Proof}

\medskip

\section{ $\mathcal I_{2,0}$ and the Birkhoff Identity}

\
\indent In this section we prove a special case of our main result that if  $\mathbf{A} \in \mathcal I_{2,0}$, then the bisemigroup $\mathbf{A}_{mj}$ satisfies (BR).   This result will play a crucial role in the proof of the main Theorem in the next section.

\begin{lema}\label{general_properties_BS_I20}
	Let $\mathbf A \in \mathcal I_{2,0}$. Then 
	\begin{enumerate}[{\rm (1)}]
		\item $((x \to (0 \to y)) \to z) \to u \approx (0 \to x) \to ((0 \to y') \to (z \to u))$, \label{260615_04}  
		\item $((0 \to x) \to y) \to z \approx (x \to y) \to (y \to z)$, \label{120919_01} 
		\item $((x \to y) \to z) \to (z \to u) \approx (0 \to x) \to ((0 \to y') \to (z \to u))$,  \label{120919_02} 
		\item $((x \to y) \to z) \to x' \approx (y \to z) \to x'$,  \label{120919_03} 
		\item $[x \to [(0 \to y') \to (z \to u)]] \to [(0 \to y) \to [(z \to u) \to (0 \to (x \to y))']] \approx (z \to u) \to (0 \to (x \to y))'$, \label{120919_04} 
		\item $(x \to (0 \to y)) \to z \approx (z \to (x \to y)) \to z$, \label{120919_05} 
		\item $[0 \to (x \to (y \to z))] \to u \approx  (0 \to x') \to ((0 \to (y \to z)) \to u),$ \label{070715_04}  
		\item $(x \to y) \to ((0 \to y) \to z) \approx (x \to y) \to z$, \label{260615_07} 
		\item $(x \to y) \to ((z \to y) \to (u \to z)) \approx  (x \to y) \to (u \to z),$ \label{070715_06} 
		\item $[(0 \to y) \to z'] \to u \approx (y \to z') \to (z' \to u)$, \label{120919_06} 
		\item $(0 \to x) \to ((y \to x) \to z) \approx (y \to x) \to z$, \label{260615_11}  
		\item $x' \to (0 \to (y \to z))' \approx x' \to (x \to (y \to z))'$, \label{120919_07} 
		\item $0 \to (x \to (y \to z)) \approx 0 \to ((x' \to y) \to z)$, \label{130919_01}
		\item $0 \to [x \to ((y \to z) \to u)] \approx 0 \to [((x \to y) \to z) \to u]$, \label{130919_02}  
		\item $(x \to y') \to [y' \to (0 \to x')'] \approx y' \to (0 \to x')'$,  \label{130919_03} 
		\item $[x \to (x' \to y)']' \approx x' \to (0 \to y')'$.  \label{130919_04} 
	\end{enumerate}
\end{lema}


\begin{Proof} Let $a,b,c,d \in A$.  

\begin{itemize}
		
		\item[(\ref{260615_04})]
		
		\noindent $(0 \to a) \to ((0 \to b') \to (c \to d)) $
		$\overset{   \ref{general_properties} (\ref{cuasiConmutativeOfImplic}) 
		}{=}  (a' \to 0') \to ((0 \to b') \to (c \to d))$
		$\overset{   \ref{general_properties}     
			(\ref{260615_01})   
		}{=}  [(0 \to a') \to (0 \to b')] \to (c \to d) $
		$\overset{   (\ref{311014_06}) \text{ and	}(\ref{071114_04})  \text{ of }  \ref{general_properties} 
		}{=}  [0 \to (a' \to b')] \to (c \to d) $
		$\overset{  (\ref{071114_04}) \text{ and } (\ref{071114_01})  \text{ of } \ref{general_properties}  
		}{=}  [0 \to (a \to b)'] \to (c \to d) $
		$\overset{   \ref{general_properties}	(\ref{cuasiConmutativeOfImplic}) 
		}{=}  [(a \to b) \to 0'] \to (c \to d) $
		$\overset{   \ref{general_properties}  (\ref{260615_01}) 
		}{=}  [[0 \to (a \to b)] \to c] \to d $
		$\overset{   \ref{general_properties}       (\ref{071114_04})
		}{=}  [[a \to (0 \to b)] \to c] \to d $.
		
		\item[(\ref{120919_01})]
		
		\noindent $((0 \to a) \to b) \to c $
		$\overset{  	\ref{general_properties} (\ref{260615_01}) 
		}{=}  (a \to 0') \to (b \to c) $
		$\overset{   \ref{general_properties}	(\ref{cuasiConmutativeOfImplic}) 
		}{=}   (0 \to a') \to (b \to c) $
		$\overset{   \ref{general_properties}	(\ref{250615_04}). 
		}{=}  (a \to b) \to (b \to c) $.
		
		\item[(\ref{120919_02})]
		\noindent $((a \to b) \to c) \to (c \to d) $
		$\overset{   \ref{general_properties} (\ref{250615_06}) 
		}{=}  (0 \to a) \to ((b \to c) \to (c \to d)) $
		$\overset{   \ref{general_properties}	(\ref{250615_04})
		}{=}  (0 \to a) \to ((0 \to b') \to (c \to d)) $.
		
		\item[(\ref{120919_03})]
		
		\noindent $((a \to b) \to c) \to a' $
		$\overset{  (I) 
		}{=}  [(a'' \to (a \to b)) \to (c \to a')']' $
		$\overset{  
		}{=}   [(a \to (a \to b)) \to (c \to a')']' $
		$\overset{   \ref{general_properties2}	(\ref{031114_04}) 
		}{=}  [(a \to b) \to (c \to a')']' $
		$\overset{  
		}{=}  [(a'' \to b) \to (c \to a')']' $
		$\overset{  (I) 
		}{=}  (b \to c) \to a' $.
		
		\item[(\ref{120919_04})]
		
		$$
		\begin{array}{cll}
		& (c \to d) \to (0 \to (a \to b))' & \mbox{} \\
		= & [[(0 \to (a \to b)) \to c] \to d] \to (0 \to (a \to b))'   &  \mbox{} \\
		& {\color{red}\qquad \qquad}\mbox{by (\ref{120919_03}) with } x:= 0 \to (a \to b), y:= c, z:= d & \mbox{} \\
		= & [[(a \to (0 \to b)) \to c] \to d] \to (0 \to (a \to b))' &  \\
		& \qquad \qquad \mbox{by Lemma \ref{general_properties2} (\ref{071114_04})} \\
		= & [(0 \to a) \to ((0 \to b') \to (c \to d))] \to (0 \to (a \to b))' &  \\
		& \qquad \qquad \mbox{by (\ref{260615_04})} \\
		= & [a \to ((0 \to b') \to (c \to d))] \to [((0 \to b') \to (c \to d)) \to (0 \to (a \to b))'] &  \\
		& {\color{red}\qquad \qquad}\mbox{by (\ref{120919_01}) with } x:= a, y:= (0 \to b') \to (c \to d), z:= (0 \to (a \to b))' & \mbox{} \\
		= & [a \to ((0 \to b') \to (c \to d))] \to [(b' \to (c \to d)) \to ((c \to d) \to (0 \to (a \to b))')] & \mbox{} \\
		& {\color{red}\qquad \qquad}\mbox{by (\ref{120919_01}) with } x:= b', y:= c \to d, z:= (0 \to (a \to b))' & \mbox{} \\
		= & [a \to ((0 \to b') \to (c \to d))] \to \{(0 \to b) \to [(0 \to 0') \to [(c \to d) \to (0 \to (a \to b))']]\} & \mbox{} \\
		& {\color{red}\qquad \qquad}\mbox{by (\ref{120919_02}) with } x:= b, y:= 0, z:= c \to d, u:= (0 \to (a \to b))' & \mbox{} \\
		= & [a \to ((0 \to b') \to (c \to d))] \to \{(0 \to b) \to [(0'' \to 0') \to [(c \to d) \to (0 \to (a \to b))']]\} & \mbox{} \\
		= & [a \to ((0 \to b') \to (c \to d))] \to \{(0 \to b) \to [0' \to [(c \to d) \to (0 \to (a \to b))']]\} &  \\
		& \qquad \qquad \mbox{by Lemma \ref{general_properties_equiv} (\ref{LeftImplicationwithtilde})} \\
		= & [a \to ((0 \to b') \to (c \to d))] \to \{(0 \to b) \to [(c \to d) \to (0 \to (a \to b))']\} &  \\
		& \qquad \qquad \mbox{by Lemma \ref{general_properties_equiv} (\ref{TXX})}
		\end{array}
		$$

		\item[(\ref{120919_05})]
		\noindent $(a \to (0 \to b)) \to c $
		$\overset{  	\ref{general_properties2} (\ref{071114_04}) 
		}{=}  (0 \to (a \to b)) \to c $
		$\overset{   \ref{general_properties2}	(\ref{291014_10}) 
		}{=}  (c \to (a \to b)) \to c $.
		
		\item[(\ref{070715_04})]

		\noindent $(0 \to a') \to ((0 \to (b \to c)) \to d) $
		$\overset{  \ref{general_properties} (\ref{cuasiConmutativeOfImplic}) 
		}{=}  (a \to 0') \to ((0 \to (b \to c)) \to d) $
		$\overset{   \ref{general_properties}    (\ref{260615_01}) 
		}{=}  [(0 \to a) \to (0 \to (b \to c))] \to d $
		$\overset{   (\ref{311014_06}) \text{ and } (\ref{071114_04}) \text{ of }  \ref{general_properties}  
		}{=}   [0 \to (a \to (b \to c))] \to d $.
		
		\item[(\ref{260615_07})]

		\noindent $		(a \to b) \to ((0 \to b) \to c) $
		$\overset{  (I) 
		}{=}  [[((0 \to b) \to c)' \to a] \to [b \to ((0 \to b) \to c)]']' $
		$\overset{   \ref{general_properties}   (\ref{250615_22}) 
		}{=}  [[((0 \to b) \to c)' \to a] \to [b \to c]']' $
		$\overset{   \ref{general_properties}   (\ref{291014_10}) 
		}{=}  [[((c \to b) \to c)' \to a] \to [b \to c]']' $
		$\overset{   \ref{general_properties} (\ref{291014_09}) 
	}{=}  [[(c \to (b \to c)')'' \to a] \to [b \to c]']' $
$\overset{  (I_{2,0}) 
}{=}  [[(c \to (b \to c)') \to a] \to [b \to c]']' $
$\overset{ \text{ (I) and (I$_{2,0}$)} 
}{=}  [[b \to c] \to (c \to (b \to c)')] \to [a \to [b \to c]']' $
$\overset{   \ref{general_properties} (\ref{281114_01}) 
}{=}  [c \to (b \to c)'] \to [a \to [b \to c]']' $
$\overset{   \ref{general_properties} (\ref{071114_05}) 
}{=}  [(b \to c) \to c'] \to [a \to [b \to c]']' $
$\overset{\text{  (I$_{2,0}$) and (I) }
}{=}  [(c' \to a) \to [b \to c]']' $
$\overset{  (I)
}{=}  (a \to b) \to c $.

		\item[(\ref{070715_06})]
		
		\noindent $(a \to b) \to ((c \to b) \to (d \to c)) $
		$\overset{  (I) 
		}{=}  (a \to b) \to [[(d \to c)' \to c] \to [b \to (d \to c)]']' $
		$\overset{   \ref{general_properties}   (\ref{250615_01}) 
		}{=}  (a \to b) \to [(d \to c) \to [b \to (d \to c)]']' $
		$\overset{   \ref{general_properties}   (\ref{291014_09}) 
		}{=}  (a \to b) \to [[(d \to c) \to b] \to (d \to c)] $
		$\overset{   \ref{general_properties}   (\ref{291014_10}) 
		}{=}  (a \to b) \to [[0 \to b] \to (d \to c)] $
		$\overset{  (\ref{260615_07})
		}{=}  (a \to b) \to (d \to c) $.
		
		\item[(\ref{120919_06})]
		
		\noindent $[(0 \to b) \to (c \to 0)] \to d $
		$\overset{   \ref{general_properties}    (\ref{260615_05}) 
		}{=}  (b \to 0') \to ((c \to 0) \to d) $
		$\overset{   \ref{general_properties}	(\ref{cuasiConmutativeOfImplic}) 
		}{=}  (0 \to b') \to ((c \to 0) \to d) $
		$\overset{   \ref{general_properties}	(\ref{250615_04})
		}{=}  (b \to c') \to (c' \to d) $.
		
		\item[(\ref{260615_11})]
		

\noindent $		(0 \to a) \to ((b \to a) \to c) $
$\overset{  \text{	\ref{general_properties} (\ref{cuasiConmutativeOfImplic})  and (I$_{2,0}$) } 
}{=}  (a' \to 0') \to ((b \to a) \to c) $
$\overset{   \ref{general_properties}   (\ref{260615_01}) 
}{=}  [(0 \to a') \to (b \to a)]  \to c $
$\overset{   \ref{general_properties}   (\ref{080415_01}) 
}{=}  [(a \to a') \to (b \to a)]  \to c $
$\overset{ \text{  \ref{general_properties_equiv}	(\ref{LeftImplicationwithtilde})  and (I$_{2,0}$) }
}{=}  [a' \to (b \to a)]  \to c $
$\overset{ \text{ \ref{general_properties} (\ref{281114_01})  and (I$_{2,0}$)
}}{=}  (b \to a)  \to c. $

		\item[(\ref{120919_07})]
		
		$$
		\begin{array}{cll}
		&  a' \to (a \to (b \to c))'  & \mbox{} \\
		= &  [(a \to (b \to c)) \to a]' \to [0 \to (a \to (b \to c))]' & \mbox{} \\
		& \qquad \qquad \mbox{by Lemma \ref{general_properties} (\ref{271114_02}) with } x:= a \to (b \to c), y:= a & \mbox{}   \\
		= & [(b \to (0 \to c)) \to a]' \to [0 \to (a \to (b \to c))]' & \mbox{} \\
		& \qquad \qquad \mbox{ by (\ref{120919_05})}  & \mbox{}   \\
		= & [(0 \to b) \to ((0 \to c') \to a')] \to [0 \to (a \to (b \to c))]'  &  \\
		& \qquad \qquad \mbox{by (\ref{260615_04}) with } x:= b, y:= c, z:= a, u:= 0 & \mbox{} \\
		= & [(0 \to b) \to ((0 \to c') \to a')] \to [(0 \to a') \to (0 \to (b \to c))'] &  \\
		& \qquad \qquad  \mbox{by (\ref{070715_04}) with } u:= 0 & \mbox{} \\
		= & (b \to ((0 \to c') \to a')) \to [((0 \to c') \to a') \to [(0 \to a') \to (0 \to (b \to c))']] & \mbox{} \\
		& \qquad \qquad \mbox{by (\ref{120919_01}) with } x:= b, y:= (0 \to c') \to a', z:= (0 \to a') \to (0 \to (b \to c))' & \mbox{} \\
		= & (b \to ((0 \to c') \to a')) \to [((0 \to c') \to a') \to (0 \to (b \to c))'] & \mbox{} \\
		& \qquad \qquad \mbox{by (\ref{070715_06}) with } x:= 0 \to c', y:= a', z:= 0, u:= 0 \to (b \to c) & \mbox{}  \\
		= & (b \to ((0 \to c') \to a')) \to [(c' \to (a \to 0)) \to [(a \to 0) \to  (0 \to (b \to c))']] & \mbox{} \\
		& \qquad \qquad \mbox{by (\ref{120919_06}) with } y:= c', z:= a, u:= (0 \to (b \to c))' & \mbox{} \\
		= & (b \to ((0 \to c') \to a')) \to [((c \to 0) \to a') \to [a' \to  (0 \to (b \to c))']] & \mbox{} \\
		= & (b \to ((0 \to c') \to a')) \to [(0 \to c) \to [(0 \to 0') \to (a' \to (0 \to (b \to c))')]] & \mbox{} \\
		& \qquad \qquad \mbox{by (\ref{120919_02}) with } x:= c, y:= 0, z:= a', u:= (0 \to (b \to c))' & \mbox{} \\
		= & (b \to ((0 \to c') \to a')) \to [(0 \to c) \to [(0'' \to 0') \to (a' \to (0 \to (b \to c))')]] & \mbox{} \\
		= & (b \to ((0 \to c') \to a')) \to [(0 \to c) \to [0' \to (a' \to (0 \to (b \to c))')]] & \mbox{ } \\
		& \qquad \qquad \mbox{by Lemma \ref{general_properties_equiv} (\ref{LeftImplicationwithtilde})}  & \mbox{}   \\
		= & (b \to ((0 \to c') \to a')) \to [(0 \to c) \to (a' \to (0 \to (b \to c))')] & \mbox{} \\
		& \qquad \qquad \mbox{ by Lemma \ref{general_properties_equiv} (\ref{TXX})}  & \mbox{}   \\
		= & (a \to 0) \to (0 \to (b \to c))' & \mbox{} \\
		& \qquad \qquad \mbox{by (\ref{120919_04}) with } x:= b, y:= c, z:= a, u:= 0 & \mbox{} \\
		= & a' \to (0 \to (b \to c))'. & \mbox{} 
		\end{array}
		$$
		
		\item[(\ref{130919_01})]
		
		\noindent $0 \to (a \to (b \to c)) $
		$\overset{  	\ref{general_properties} (\ref{140715_20}) 
		}{=}  0 \to [(a \to b')' \to c] $
		$\overset{\text{  \ref{general_properties2}	(\ref{071114_04}) and (\ref{311014_06}) }
		}{=}  [0 \to (a \to b')'] \to (0 \to c) $
		$\overset{  \ref{general_properties2} (\ref{191114_05}) 
		}{=}  [0 \to (a' \to b)] \to (0 \to c) $
		$\overset{ \text{  \ref{general_properties2}	(\ref{311014_06})  and  (\ref{071114_04}) }
		}{=}  0 \to ((a' \to b) \to c) $.
		
		\item[(\ref{130919_02})]
		
		\noindent $0 \to [a \to ((b \to c) \to d)]  $
		$\overset{ 	(\ref{130919_01}) 
		}{=}  0 \to [(a' \to (b \to c)) \to d]  $
		$\overset{ \text{  \ref{general_properties2}	(\ref{071114_04}) and (\ref{311014_06}) }
		}{=}  [0 \to (a' \to (b \to c))] \to (0 \to d) $
		$\overset{  (\ref{130919_01}) 
		}{=}  [0 \to ((a'' \to b) \to c)] \to (0 \to d) $
		$\overset{  
		}{=}  [0 \to ((a \to b) \to c)] \to (0 \to d) $
		$\overset{\text{   \ref{general_properties2}	(\ref{311014_06})  and  (\ref{071114_04})}
		}{=}  0 \to [((a \to b) \to c) \to d] $.
		
		\item[(\ref{130919_03})]
		
		\noindent $(a \to b') \to [b' \to (0 \to a')'] $
		$\overset{ 	 \ref{general_properties2} (\ref{250615_04}) 
		}{=}  (0 \to a') \to [b' \to (0 \to a')'] $
		$\overset{   \ref{general_properties2} (\ref{281114_01}) 
		}{=}  b' \to (0 \to a')' $.
		
		\item[(\ref{130919_04})]

		$$
		\begin{array}{lcll}
		[a \to (a' \to b)']' & = & [(a' \to 0) \to (a' \to b)']' & \mbox{} \\
		& = & [[(a' \to b)'' \to a'] \to (0 \to (a' \to b)')']'' &  \\
		&  & \qquad \qquad \mbox{by (I)} &  \\
		& = & [(a' \to b) \to a'] \to (0 \to (a' \to b)')' & \mbox{} \\
		& = & [(0 \to b) \to a'] \to (0 \to (a' \to b)')' &  \\
		&  & \qquad \qquad \mbox{by Lemma \ref{general_properties2} (\ref{291014_10})} &  \\
		& = & [(0 \to b) \to a'] \to [0 \to (a \to (0 \to b)')]' &  \\
		&  & \qquad \qquad \mbox{by (\ref{130919_02}) with } x:=a, y:= 0, z:= b, u:= 0 & \mbox{} \\
		& = & [(0 \to b) \to a'] \to [a \to (0 \to (0 \to b)')]' &  \\
		&  & \qquad \qquad \mbox{by Lemma \ref{general_properties2} (\ref{071114_04})} &  \\
		& = & [(0 \to b) \to a'] \to [a \to (0 \to b')]' &  \\
		&  & \qquad \qquad \mbox{by Lemma \ref{general_properties2} (\ref{031114_07})} &  \\
		& = & [(0 \to b) \to a'] \to [a'' \to (0 \to b')]' &   \mbox{ }\\ 
		& = & [(0 \to b) \to a'] \to [(0 \to a') \to (0 \to b')'] &  \\
		&  & \qquad \qquad \mbox{by Lemma \ref{general_properties2} (\ref{031114_06})} &  \\
		& = & [(0 \to b) \to a']  \to (0 \to b')' & \mbox{} \\
		&  & \qquad \qquad \mbox{by (\ref{070715_06}) with } x:= 0 \to b, y:= a', z:= 0, u:= 0 \to b' & \mbox{} \\
		& = & (b \to a') \to [a' \to (0 \to b')'] &  \\
		&  & \qquad \qquad \mbox{by (\ref{120919_06}) with } y:= b, z:= a, u:= (0 \to b')' & \mbox{} \\
		& = & a' \to (0 \to b')' &  \\
		&  & \qquad \qquad \mbox{by (\ref{130919_03}) with } x := b, y:= a. &  
		\end{array}
		$$		
	\end{itemize}		
\end{Proof}

The proof of our main result (Theorem \ref{Main}), given in the next section, depends on the following theorem. 
 
\begin{Theorem} \label{theo_BS_I20}
	Let $\mathbf A \in \mathcal I_{2,0}$.  Then $\mathbf A_{mj}$ satisfies the Birkhoff identity. 
\end{Theorem}

\begin{Proof}
Let $a,b \in A$. Then


$$
\begin{array}{lcll}
a \wedge (a \vee b) 
& = & (a \to (a' \to b)')' & \mbox{by definition of } \vee  \mbox{ and } \wedge \\
& = & a' \to (0 \to b')' & \mbox{by Lemma \ref{general_properties_BS_I20} (\ref{130919_04})} \\
& = & a' \to (a \to b')' & \mbox{by Lemma \ref{general_properties_BS_I20} (\ref{120919_07}) with } z:= 0 \\
& = & [a' \to (a \to b')']'' & \mbox{} \\
& = & (a' \wedge (a \wedge b)')' & \mbox{by definition of } \wedge \\
& = & a \vee (a \wedge b) & \mbox{by definition of } \wedge. 
\end{array}
$$	
\end{Proof}

\section{Main theorem}

\
\indent In this section, the main theorem of this paper is proved.  For this we need one more crucial result proved in \cite{cornejo2016derived}.

\begin{Theorem}[Transfer Theorem] {\rm \cite{cornejo2016derived}}  \label{Theo_identities_I20_to_I}
	
	Let $t_i(\overline x), i= 1, \cdots, 6$.
	be terms, where $\overline x$ denotes the sequence $\langle x_1, \cdots x_n \rangle$, $x_i$ being varaibles.  Let $\mathcal V$ be a subvariety of $\mathcal I$.  \\ 
	If 
	$$\mathcal V  \cap \mathcal I_{2,0} \models (t_1(\overline x) \to t_2(\overline x)) \to t_3(\overline x) \approx (t_4(\overline x) \to t_5(\overline x)) \to t_6(\overline x),$$ 
	then 
	$$\mathcal V\models (t_1(\overline x) \to t_2(\overline x)) \to t_3(\overline x) \approx (t_4(\overline x) \to t_5(\overline x)) \to t_6(\overline x).$$
\end{Theorem}

We are now ready to present the main result of this paper (i.e., Theorem \ref{main} of Introduction).

\begin{Theorem} \label{Main}
	Let $\mathbf A \in \mathcal I$ then $\mathbf A_{mj}$ satisfies the Birkhoff's identity:
	\[{\rm (BR)} \quad x \land (x \lor y) \approx x \lor (x \land y).\]  
\end{Theorem}

\begin{Proof}
Apply Theorem \ref{theo_BS_I20} and Theorem \ref{Theo_identities_I20_to_I}.	
\end{Proof}

\ \\
\indent  Recall that $\mathcal{S}= \mathcal{I}_{2,0} \cap \mathcal{MC}$. 
The following corollary, which characterizes the implication zroupoids $\mathbf A$ for which  $\mathbf A_{mj}$ is a Birkhoff system, is an improvement on \cite[Theorem 7.3]{cornejo2015implication}. 

\begin{Corollary} \label{Theo_BS_iff}
	Let $\mathbf A \in \mathcal I$.  Then the algebra $\mathbf A_{mj}$ 
	is a Birkhoff system if and only if $\mathbf A \in \mathcal{S}$. 
\end{Corollary}

\begin{Proof}
Let  $\mathbf A_{mj}$  
be a Birkhoff system, then $\mathbf A$ satisfies (MC) and the identity: $x \wedge x \approx x$, which, in view of Lemma \ref{general_properties_equiv}, implies that $\mathbf A \models x \approx x''$.  Hence $\mathbf A \in \mathcal S$.
	For the converse, assume that $\mathbf A \in \mathcal{S}$. Then $\mathbf A$ satisfies the identities $x \wedge y \approx y \wedge x$ and $x \approx x''$. We know also by Theorem \ref{Theo_Assoc} that the operations $\land$ and $\vee$ are associative. 
	Hence, from Theorem \ref{Main} (or Theorem \ref{theo_BS_I20}) we conclude that  $\mathbf A_{mj}$  
is a Birkhoff system.
\end{Proof}
\ \\ 

\indent Recall that implication zroupoids that satisfy the associative identity:  $$(A) \ \ \ x \to (y \to z) \approx (x \to y) \to z$$were called \emph {implication semigroups} in \cite{Implsem}.  Let  $\mathcal{IS}$ denote the variety of implication semigroups.  The following special case of our main result may be of interest to the researchers in semigroup theory.

\begin{Corollary} \label{corA}
Let $\mathbf A \in \mathcal{IS}$.  Then the algebra $\mathbf A_{mj}$ satisfies the Birkhoff identity.
\end{Corollary}

In the next section we will improve the above corollary.






\section{Concluding Remarks}
\

As mentioned in the introduction, it follows from the main result, Theorem \ref{Main}, that the algebras $\mathbf A \in \mathcal{I \setminus S}$ such that $\mathbf A_{mj}$ are bisemigroups that are not bisemilattices and satisfy the Birkhoff's identity, which suggests naturally  a  generalization of Birkhoff systems, which we will call ``Birkhoff bisemigroups''.  To the best of our knowledge, the algebras defined in the following definition seem to be new.


\begin{definition}
A bisemigroup $\mathbf{A}$ is a Birkhoff bisemigroup if $\mathbf{A} \models (BR)$.
\end{definition}

 In this new terminology, we can now recast our main theorem as: If $\mathbf{A} \in \mathcal I$, 
then $\mathbf{A}_{mj}$ is a Birkhoff bisemigroup.

Thus, the class of algebras $\mathbf{A_{mj}}$, where $\mathbf{A} \in \mathcal I$ 
provide a large class of examples of Birkhoff bisemigroups. 

Another class of examples of Birkhoff bisemigroups arise from semigroups themselves as follows: Let $\mathbf A=\langle A, \land \rangle$ be a semigroup. Then the algebra $\langle A, \land, \land \rangle$  is clearly a bisemigroup and satisfies the Birkhoff identity trivially as the two binary operations are the same.     
Let us call such a bisemigroup arising from a semigroup ``essentially a semigroup''.  We shall now improve and clarify Corollary \ref{corA}.  For this we need the following lemma:

\begin{lema}
Let $\mathbf A \in \mathcal{IS}$.  Then $\mathbf A$ holds:
	\begin{enumerate}[{\rm (1)}]
		\item $0 \to 0'\approx 0$, \label{130120_01}
		\item $0 \to x'\approx x'$, \label{130120_02}
		\item $0'\approx 0$, \label{130120_03}
		\item $x \lor y \approx x \land y$. \label{130120_04}
	\end{enumerate}
\end{lema}

\begin{Proof} Let $a,b \in A$.  
	
	\begin{itemize}
		
		\item[(\ref{130120_01})]
\noindent $	0	$
$\overset{  
}{=}  0'' $
$\overset{  
}{=}  (0 \to 0)\to 0 $
$\overset{  (A) 
}{=}  0 \to (0 \to 0) $
$\overset{  
}{=}  0 \to 0'. $
		
		\item[(\ref{130120_02})]
\noindent $a'	$
$\overset{  
}{=}  a \to 0 $
$\overset{  (\ref{130120_01}) 
}{=}  a \to (0 \to 0') $
$\overset{  (A) 
}{=}  (a \to 0) \to 0' $
$\overset{  (I) 
}{=}  [(0'' \to a) \to (0 \to 0')']' $
$\overset{  
}{=}  [(0 \to a) \to (0 \to 0')']' $
$\overset{  (\ref{130120_01}) 
}{=}  [(0 \to a) \to 0'] \to 0 $
$\overset{  (A) 
}{=}  (0 \to a) \to (0' \to 0) $
${=}  (0 \to a) \to 0'' 
{=}  (0 \to a) \to 0 $
$\overset{  (A) 
}{=}  0 \to (a \to 0) $
$\overset{  
}{=}  0 \to a' $.

		\item[(\ref{130120_03})]
\noindent $0' $
$\overset{  (\ref{130120_02}) 
}{=}  0 \to 0' $
$\overset{  (\ref{130120_01}) 
}{=}  0 $.	
		
		\item[(\ref{130120_04})]
\noindent $a \lor b $
$\overset{  \mbox{ def of } \lor \mbox{ and } (A) 
}{=}  a \to (0 \to (b \to (0 \to (0 \to (0 \to 0))))) $
$\overset{  (\ref{130120_03}) 
}{=}  a \to (0 \to (b \to 0)) $
$\overset{  (\ref{130120_02}) 
}{=}  a\to (b \to 0) $
$\overset{  (\ref{130120_03}) 
}{=}  a \to (b \to (0 \to 0)) $
$\overset{  (A) 
}{=}  a \to [(b \to 0) \to 0] $
$\overset{  (A) 
}{=}  [a \to (b \to 0)] \to 0 $
$\overset{  \mbox{ def of } \land 
}{=}  a \land b $.			
	\end{itemize}
\end{Proof}



In view of the preceding lemma, Corollary \ref{corA} can be improved to the following.
\begin{Corollary} 
Let $\mathbf A \in \mathcal{IS}$.  Then the algebra $\mathbf A_{mj}$ is essentially a semigroup.
\end{Corollary}

We conclude this paper with the following open problem. \\

{\bf PROBLEM:} Investigate Birkhoff bisemigroups; in particular, describe the structure of the lattice of subvarieties of the variety of Birkhoff bisemigroups.

\section*{Acknowledgements} 

 The first author wants to thank the institutional support of CONICET (Consejo Nacional de Investigaciones Cient\'ificas y T\'ecnicas) and Universidad Nacional del Sur..  The authors also wish to acknowledge that [Mc] was a useful tool during the research phase of this paper.\\

\noindent {\bf Compliance with Ethical Standards:}\\ 

{\bf Conflict of Interest}: The first author declares that he has no conflict of interest. The second author declares that he has no conflict of interest.\\

{\bf Ethical approval}: This article does not contain any studies with human participants or animals performed by any of the authors.\\


\noindent {\bf Funding:}  
	The work of Juan M. Cornejo was supported by CONICET (Consejo Nacional de Investigaciones Cientificas y Tecnicas) and Universidad Nacional del Sur.

\vskip 1.5cm

\noindent {\sc Juan M. Cornejo}\\
Departamento de Matem\'atica\\
Universidad Nacional del Sur\\
Alem 1253, Bah\'ia Blanca, Argentina\\
INMABB - CONICET\\

\noindent jmcornejo@uns.edu.ar

\vskip 1.4cm

\noindent {\sc Hanamantagouda P. Sankappanavar}\\
Department of Mathematics\\
State University of New York\\
New Paltz, New York 12561\\
U.S.A.\\

\noindent sankapph@newpaltz.edu

\end{document}